\documentclass{amsart}

\theoremstyle{definition}

\theoremstyle{remark}

\reversemarginpar

\begin{document}

\title[A logarithmic family]{The integrals in Gradshteyn and Rhyzik. Part 1: \\
a family of logarithmic integrals}

\author{Victor H. Moll}
\address{Department of Mathematics,
Tulane University, New Orleans, LA 70118}
\email{vhm@math.tulane.edu}

\subjclass{Primary 33}

\date{\today}


\begin{abstract}
We present the evaluation of  a family of logarithmic integrals.  This 
provides a unified proof of several formulas
in the classical table of integrals by I. S. 
Gradshteyn and I. M. Rhyzik. 
\end{abstract}

\maketitle

\newcommand{\nn}{\nonumber}
\newcommand{\ba}{\begin{eqnarray}}
\newcommand{\ea}{\end{eqnarray}}
\newcommand{\ift}{\int_{0}^{\infty}}
\newcommand{\ione}{\int_{0}^{1}}
\newcommand{\ifft}{\int_{- \infty}^{\infty}}
\newcommand{\no}{\noindent}
\newcommand{\Ftwo}{{{_{2}F_{1}}}}
\newcommand{\realpart}{\mathop{\rm Re}\nolimits}
\newcommand{\imagpart}{\mathop{\rm Im}\nolimits}

\newtheorem{Definition}{\bf Definition}[section]
\newtheorem{Thm}[Definition]{\bf Theorem} 
\newtheorem{Example}[Definition]{\bf Example} 
\newtheorem{Lem}[Definition]{\bf Lemma} 
\newtheorem{Note}[Definition]{\bf Note} 
\newtheorem{Cor}[Definition]{\bf Corollary} 
\newtheorem{Prop}[Definition]{\bf Proposition} 
\newtheorem{Problem}[Definition]{\bf Problem} 
\numberwithin{equation}{section}

\section{Introduction} \label{intro} 
\setcounter{equation}{0}

The values of many definite integrals have been compiled in the classical 
{\em Table of Integrals, Series and Products} by I. S. Gradshteyn and 
I. M. Rhyzik \cite{gr}. The table is organized like a phonebook: integrals 
that {\em look} similar are place close together. For 
example, $\mathbf{4.229.4}$ gives
\begin{equation}
\ione \ln \left( \ln \frac{1}{x} \right) 
\left( \ln \frac{1}{x} \right)^{u-1} \, dx = \psi(\mu) \Gamma(\mu), \label{one}
\end{equation}
\noindent
for $\realpart{\mu} > 0$, and 
$\mathbf{4.229.7}$ states that 
\begin{equation}
\int_{\pi/4}^{\pi/2} \ln \ln \tan x \, dx = 
\frac{\pi}{2} \ln \left\{ \frac{\Gamma \left( \tfrac{3}{4} \right)} 
{\Gamma \left( \tfrac{1}{4} \right)} \, \sqrt{2 \pi} \, \right\}. \label{two}
\end{equation}
\noindent
In spite of a large amount of work in the development of this table, the 
latest version of \cite{gr} still contains some typos. For
example, the exponent $u$ in (\ref{one}) should be $\mu$. A list of 
errors and typos can be found in 
\noindent
\begin{verbatim}          http://www.mathtable.com/errata/gr6_errata.pdf\end{verbatim}

\medskip

The fact that two integrals are close in the table is not a reflection of 
the difficulty involved in their evaluation. Indeed, the formula (\ref{one}) 
can be established by the change of variables $v = - \ln x$ followed by 
differentiating the classical 
gamma function 
\begin{equation}
\Gamma(\mu) := \ift t^{\mu-1} e^{-t} \, dt, \quad \realpart{\mu} > 0,
\end{equation}
\noindent
with respect to the parameter $\mu$. The 
function $\psi(\mu)$ in (\ref{one})
is simply the logarithmic derivative of $\Gamma(\mu)$ and the formula has
been checked. The situation is quite different for (\ref{two}). This formula 
is the subject of the lovely paper \cite{vardi1} in which the author uses 
Analytic Number Theory to check (\ref{two}). The ingredients of the proof are 
quite formidable: the author shows that 
\begin{equation}
\int_{\pi/4}^{\pi/2} \ln \ln \tan x \, dx = 
\frac{d}{ds} \Gamma(s) L(s) \text{ at } s =1,
\label{valueL1}
\end{equation}
\noindent
where 
\begin{equation}
L(s) = 1 - \frac{1}{3^{s}} + \frac{1}{5^{s}} - \frac{1}{7^{s}} + \cdots 
\end{equation}
\noindent
is the Dirichlet L-function. The computation of (\ref{valueL1}) is done 
in terms of the Hurwitz zeta function
\begin{equation}
\zeta(q,s) = \sum_{n=0}^{\infty} \frac{1}{(n+q)^{s}},
\end{equation}
\noindent
defined for $0 < q < 1$ and $\realpart{s} > 1$. The function $\zeta(q,s)$ 
can be analytically continued to the whole plane with only a simple pole
at $s=1$ using the integral representation
\begin{equation}
\zeta(q,s) = \frac{1}{\Gamma(s)} \ift \frac{e^{-qt} t^{s-1}}{1-e^{-t}} \, dt.
\end{equation}
\noindent
The relation with the $L$-functions is provided by
employing
\begin{equation}
L(s) = 2^{-2s} \left( \zeta(s, \tfrac{1}{4}) - \zeta(s, \tfrac{3}{4}) \right). 
\end{equation}
\noindent
The functional equation 
\begin{equation}
L(1-s) = \left( \frac{2}{\pi} \right)^{s} \sin \frac{\pi s}{2} 
\Gamma(s) \, L(s),
\end{equation}
\noindent
and Lerch's identity 
\begin{equation}
\zeta'(0,a) = \log \frac{\Gamma(a)}{\sqrt{2 \pi}},
\end{equation}
\noindent
complete the evaluation.  More information about these functions can be 
found in \cite{ww}. \\

\medskip

In the introduction to \cite{irrbook} we expressed the desire to 
establish {\em all} the 
formulas in \cite{gr}. This is a {\em nearly impossible task} as was also 
noted by a (not so) favorable review given in \cite{shawyer1}. 
This is the first of a series of papers where we present some of these 
evaluations. 

We consider here the family
\begin{equation}
f_{n}(a) = \ift \frac{\ln^{n-1}x \, dx}{(x-1)(x+a)}, \text{ for } n \geq 2 
\text{ and } a>0. 
\label{fofn}
\end{equation}
\noindent 
Special examples of $f_{n}$ appear in \cite{gr}. The reader will find 
\begin{equation}
f_{2}(a) = \frac{\pi^{2} + \ln^{2}a}{2(1+a)} 
\end{equation}
as formula $\mathbf{4.232.3}$ and 
\begin{equation}
f_{3}(a) = \frac{\ln a \, ( \pi^{2} + \ln^{2}a)}{3(1+a)} 
\end{equation}
as formula $\mathbf{4.261.4}$. In later sections the persistent reader will find 
\begin{eqnarray}
f_{4}(a) & = & \frac{(\pi^{2} + \ln^{2}a)^{2}}{4(1+a)} \nonumber \\
f_{5}(a) & = & \frac{\ln a \, (\pi^{2} + \ln^{2}a) (7 \pi^{2} + 3 \ln^{2}a)}
{15(1+a)} \nonumber \\
f_{6}(a) & = & \frac{(\pi^{2} + \ln^{2}a)^{2}(3 \pi^{2} + \ln^{2}a)}{6(1+a)} 
\nonumber
\end{eqnarray}
\noindent
as $\mathbf{4.262.3}, \, \mathbf{4.263.1}$ and $\mathbf{4.264.3}$ respectively. 

These formulas suggest that
\begin{equation}
h_{n}(b) := f_{n}(a) \times (1+a)
\label{polyh}
\end{equation}
\noindent 
is a polynomial in the variable $b = \ln a$. The relatively elementary 
evaluation of $f_{n}(a)$ discussed here identifies this polynomial. 

There are several classical results that are stated without proof. The reader 
will find them in \cite{apostol1} and \cite{irrbook}. 

\section{The evaluation} \label{evaluation} 
\setcounter{equation}{0}

The expression (\ref{fofn}) for $f_{n}(a)$ can be written as
\begin{equation}
f_{n}(a) = \int_{0}^{1} \frac{\ln^{n-1}x \, dx}{(x-1)(x+a)}
+ \int_{1}^{\infty} \frac{\ln^{n-1}x \, dx}{(x-1)(x+a)},
\nonumber
\end{equation}
\noindent
and the transformation $t = 1/x$ in the second integral yields 
\begin{equation}
f_{n}(a) = \int_{0}^{1} \frac{\ln^{n-1}x \, dx}{(x-1)(x+a)}
+ (-1)^{n}\int_{0}^{1} \frac{\ln^{n-1}x \, dx}{(x-1)(1+ax)}.
\nonumber
\end{equation}
\noindent 
The partial decomposition
\begin{equation}
\frac{1}{(x-1)(x+a)} = \frac{1}{1+a} \frac{1}{x-1} - \frac{1}{1+a} \frac{1}{x+a}
\nonumber 
\end{equation}
\noindent
yields the representation
\begin{equation}
f_{n}(a) =  \frac{1 - (-1)^{n-1}}{1+a} \int_{0}^{1} \frac{\ln^{n-1}x \, dx}{x-1} - 
 \frac{1}{1+a} \int_{0}^{1} \frac{\ln^{n-1}x \, dx}{x+a}  
+ (-1)^{n-1} \frac{a}{1+a} 
 \int_{0}^{1} \frac{\ln^{n-1}x \, dx}{1 + ax}.  \nonumber 
\end{equation}

The  evaluation of these integrals require the {\em polylogarithm} function 
defined by 
\begin{equation}
\text{Li}_{m}(x):= \sum_{k=1}^{\infty} \frac{x^{k}}{k^{m}}. 
\end{equation}
\noindent
This function is sometimes denoted by $\text{PolyLog}[m,x]$. 
Detailed information 
about the polylogarithm functions appears in \cite{lewin1}. 

\begin{Prop}
For  $n \in \mathbb{N}, \, n \geq 2$ and $a > 1$ we have 
\begin{eqnarray}
\ione \frac{\ln^{n-1}x \, dx}{x-1}  & = & (-1)^{n} (n-1)! \, \zeta(n),  
\nonumber \\
\ione \frac{\ln^{n-1}x \, dx}{x+a}  & = & (-1)^{n} (n-1)! \,
\text{Li}_{n}(-1/a),   \nonumber \\
\ione \frac{\ln^{n-1}x \, dx}{1+ ax}  & = & (-1)^{n} \frac{(n-1)!}{a}  \, 
\text{Li}_{n}(-a). \nonumber 
\end{eqnarray}
\end{Prop}
\begin{proof}
Simply expand the integrand in a geometric series.
\end{proof}

\medskip

\begin{Cor}
\label{cor1}
The integral $f_{n}(a)$ is given by 
\begin{equation}
f_{n}(a) = \frac{(-1)^{n} (n-1)!}{1+a} 
\left\{ \left[ ( 1 - (-1)^{n-1} \right] \zeta(n) - \text{Li}_{n} 
\left(-\tfrac{1}{a} \right) + 
(-1)^{n-1} \text{Li}_{n}(-a) \right\}. \nonumber
\end{equation}
\end{Cor}

\bigskip

The reduction of the previous expression requires the identity
\begin{equation}
\text{Li}_{\nu}(z) =  \frac{(2 \pi)^{\nu}}{\Gamma(\nu)} 
e^{\pi i \nu/2} \zeta \left( 1 - \nu, \frac{\log(-z)}{2 \pi i } + \frac{1}{2}
\right)  - 
    e^{\pi i \nu} \text{Li}_{\nu}(-1/z). 
\end{equation}
\noindent
This transformation for the polylogarithm function appears in 

\begin{verbatim}
              http://functions.wolfram.com/10.08.17.0007.01 
\end{verbatim}

\medskip

In the special case 
$z = -a$ and $\nu = n$, with $n \in \mathbb{N}, \, n \geq 2$, 
 we obtain 
\begin{equation}
(-1)^{n-1} \text{Li}_{n}(-a) - \text{Li}_{n}(-1/a) = 
\frac{(2 \pi)^{n} }{n! \, i^{n}} B_{n} \left( \frac{\log a}{2 \pi 
i} + \frac{1}{2} \right),
\label{reduction}
\end{equation}
\noindent
where $B_{n}(z)$ is the Bernoulli polynomial of order $n$. This family of 
polynomials is defined by their exponential generating function
\begin{equation}
\frac{t e^{qt}}{e^{t} - 1} = \sum_{k=0}^{\infty} B_{k}(q) \frac{t^{k}}{k!}.
\end{equation}
\noindent
The classical identity
\begin{equation}
\zeta(1-k,q) = - \frac{1}{k} B_{k}(q), \text{ for } k \in \mathbb{N}
\end{equation}
\noindent
is used in (\ref{reduction}). Therefore the result in Corollary \ref{cor1}
can be written as:

\begin{Cor}
\label{cor2}
The integral $f_{n}(a)$ is given by 
\begin{equation}
f_{n}(a) = \frac{(-1)^{n}}{1+a} (n-1)! \left[ 1 + (-1)^{n} \right] \zeta(n)
+ \frac{(2 \pi i)^{n}}{n(1+a)} B_{n} \left( \frac{\log a }{2 \pi i} 
+ \frac{1}{2} \right). \nonumber
\end{equation}
\end{Cor}

\medskip

We now proceed to simplify this representation.  The 
Bernoulli polynomials satisfy the addition theorem 
\begin{equation}
B_{n}(x+y) = \sum_{j=0}^{n} \binom{n}{j} B_{j}(x) y^{n-j},
\end{equation}
\noindent 
and the reflection formula 
\begin{equation}
B_{n}(\tfrac{1}{2}-x) = (-1)^{n} B_{n} \left( \tfrac{1}{2} + x \right). 
\end{equation}
\noindent 
In particular $B_{n}( \tfrac{1}{2} ) = 0$ if $n$ is odd. For $n$ even, one 
has 
\begin{equation}
B_{n}( \tfrac{1}{2} ) = ( 2^{1-n}-1)B_{n},
\end{equation}
\noindent
where $B_{n}$ is the Bernoulli number $B_{n}(0)$. Thus, the last 
term in Corollary 
\ref{cor2} becomes
\begin{equation}
B_{n} \left( \frac{\log a}{2 \pi i} + \frac{1}{2} \right) = 
\sum_{j=0}^{\lfloor{n/2 \rfloor}} \binom{n}{2j} (2^{1-2j}-1)B_{2j}
\left( \frac{\log a}{2 \pi i } \right)^{n-2j}. 
\nonumber
\end{equation}

We have completed the proof of the following closed-form formula for 
$f_{n}(a)$: \\ 

\begin{Thm}
The integral $f_{n}(a)$ is given by 
\begin{eqnarray}
f_{n}(a) & = & \frac{(-1)^{n} \, (n-1)!}{1+a} \left[ 1 + (-1)^{n} \right] 
\zeta(n) + \nonumber \\
 & +  & \frac{1}{n(1+a)} 
\sum_{j=0}^{\lfloor{n/2 \rfloor}} \binom{n}{2j} 
(2^{2j}-2)(-1)^{j-1}B_{2j}  \pi^{2j} 
(\log a)^{n-2j}. 
\nonumber
\end{eqnarray}
\end{Thm}

\medskip

Observe that if $n$ is odd, the first term vanishes and there is no 
contribution of the {\em odd zeta values}. For $n$ even, the first term 
provides a rational multiple of $\pi^{n}$ in view of Euler's representation
of the even zeta values 
\begin{equation}
\zeta(2m) = \frac{(-1)^{m+1} \, (2 \pi)^{2m} B_{2m} }{2(2m)!}. 
\end{equation}

\medskip

The polynomial $h_{n}$ predicted in (\ref{polyh}) can now be read directly 
from this expression for the integral $f_{n}$. Observe that $h_{n}$ has 
positive coefficients because the Bernoulli numbers satisfy $(-1)^{j-1}B_{2j} > 0$. 

\medskip

\noindent
{\bf Note}. The change of variables $t = \ln x$ converts $h_{n}(a)$ into the
form
\begin{equation}
h_{n}(a) = \ifft \frac{t^{n-1} \, dt}{(1-e^{-t})(a+e^{t})}. 
\end{equation}
\noindent
The integrals $h_{n}(a)$  for $n=2, \, \cdots, 5$ appear in \cite{gr} as 
$\mathbf{3.419.2}, \cdots, \mathbf{3.419.6}$. The latest edition has an error in the 
expression for this last value. \\

\noindent
{\bf Conclusions}. We have provided an evaluation of the integral 
\begin{equation}
f_{n}(a) := \ift \frac{\ln^{n-1}x \, dx}{(x-1)(x+a)},
\label{fofna}
\end{equation}
\noindent
given by
\begin{eqnarray}
n(1+a)f_{n}(a) & = & (-1)^{n}n! \left[ 1 + (-1)^{n} \right] \zeta(n)
\label{formula} \\
 & + & \sum_{j=0}^{\lfloor{\tfrac{n}{2} \rfloor} } 
\binom{n}{2j} (2^{2j}-2) (-1)^{j-1}B_{2j} \pi^{2j} (\log a)^{n-2j}. 
\nonumber
\end{eqnarray}

\medskip

\noindent
{\bf Symbolic calculation}. We now describe our attempts to evaluate the 
integral $f_{n}(a)$ using Mathematica $5.2$.  For a specific value of $n$, 
Mathematica is capable of producing the result in (\ref{formula}). The 
integral is returned unevaluated if $n$ is given as a parameter. 

\bigskip

\noindent
{\bf Acknowledgments}. The author wishes to thank Luis Medina for a 
careful reading of an earlier version of the paper. The partial support of 
$\text{NSF-DMS } 0409968$ is also acknowledged.

\end{document}